\documentclass[12pt]{article}

\usepackage{amsmath}
\usepackage{amssymb}
\usepackage[mathscr]{eucal}
\usepackage{graphicx}
\usepackage{color}

\newtheorem{Theorem}{\sc Theorem}
\newtheorem{Definition}[Theorem]{\sc Definition}
\newtheorem{Proposition}[Theorem]{\sc Proposition}
\newtheorem{Lemma}[Theorem]{\sc Lemma}

\newtheorem{Example}[Theorem]{\sc Example}

\newcommand{\R}{{\if mm {\rm I}\mkern -3mu{\rm R}\else \leavevmode
\hbox{I}\kern -.17em\hbox{R} \fi}}

\newcommand{\cF}{\mbox{{${\cal F}$}}}

\newcommand{\bu}{\mbox{\boldmath{$u$}}}
\newcommand{\bv}{\mbox{\boldmath{$v$}}}
\newcommand{\bw}{\mbox{\boldmath{$w$}}}
\newcommand{\bx}{\mbox{\boldmath{$x$}}}

\newcommand{\fb}{\mbox{\boldmath{$f$}}}

\newcommand{\bsigma}{\mbox{\boldmath{$\sigma$}}}
\newcommand{\btau}{\mbox{\boldmath{$\tau$}}}
\newcommand{\bvarepsilon}{\mbox{\boldmath{$\varepsilon$}}}
\newcommand{\bnu}{\mbox{\boldmath{$\nu$}}}

\newcommand{\bzero}{\mbox{\boldmath{$0$}}}

\def\sqr#1#2{{
    \vcenter{
         \vbox{\hrule height.#2pt
               \hbox{\vrule width.#2pt height#1pt \kern#1pt
                     \vrule width.#2pt
               }
               \hrule height.#2pt
         }
    }
}}

\def\real{\mathbb{R}}

\def\lista#1
{{ \itemindent 0.0cm \labelsep .2cm \leftmargin 0.8cm \rightmargin
0.0cm \labelwidth 0.6cm \topsep 0.0mm
\parsep 0.0mm
\itemsep 0.0mm
\begin{list}{}
{ \setlength{\leftmargin}{.8cm} \setlength{\rightmargin}{0.0cm}
\setlength{\parsep}{0.0mm} \setlength{\topsep}{.0mm}
\setlength{\parskip}{.0cm} \setlength{\itemsep}{.0cm} }
{#1}\end{list}} }

\newcounter{theorem}

\begin{document}

\title{\bf On the Optimal Control of Variational-Hemivariational Inequalities}

{\author{Yi-bin Xiao$^{1}$\, and\, Mircea Sofonea$^{1,2}$\footnote{Corresponding author, E-mail :sofonea@univ-perp.fr}\\[5mm]
			{\it\small $^1$ School of Mathematical Sciences}\\
		{\it \small University of Electronic Science and Technology of China}\\
		{\it \small Chengdu, Sichuan, 611731, PR China}
		\\[6mm]	
		{\it \small $^2$ Laboratoire de Math\'ematiques et Physique}\\
		{\it \small
			University of Perpignan Via Domitia}
		\\{\it\small 52 Avenue Paul Alduy, 66860 Perpignan, France}
	}

\date{}
\maketitle
\thispagestyle{empty}

\vskip 5mm

\noindent {\bf Abstract.} \
The present paper represents a continuation of \cite{S1}. There,  a continuous dependence  result for the solution of an elliptic variational-hemivariational inequality was obtained and then used to prove the existence of optimal pairs for two associated optimal control problems.
In the current paper we  complete this study with more general results. Indeed,  we prove the continuous dependence of the solution with respect to a parameter which appears in all the data of the problem, including the set of constraints, the nonlinear operator and the two functionals which govern the variational-hemivariational inequality. This allows us to consider a general associated optimal control problem for which  we prove the existence of  optimal pairs, together with a new convergence result. 
The mathematical tools developed  in this paper are useful in the analysis and control of a large class of  boundary value problems  which, in a weak formulation, lead to  elliptic variational-hemivariational inequalities. To provide an example, we illustrate our results in the study of an inequality  which describes the equilibrium
of an elastic body in frictional contact with a foundation made of a rigid body covered by a layer of soft material.

\vskip 5mm

\noindent {\bf Keywords:} Variational-hemivariational inequality,
Clarke subdifferential, Mosco convergence, optimal pair, optimal control, elastic body, frictional contact problem.

\vskip 5mm

\noindent {\bf 2010 Mathematics Subject Classification:} \ 47J20,
47J22, 49J53, 74M15.

\vskip 15mm

\section{Introduction}\label{s1}
\setcounter{equation}0


Variational and hemivariational inequalities are widely used in the study of many nonlinear boundary value problems and have a large number of applications in Contact Mechanics and Engineering see, for instance, \cite{C,MOSBOOK,Pa,P,SX}. The theory of variational inequalities was developed in early sixty's, by using
arguments of monotonicity and convexity, including properties of the subdifferential 
of a convex function.  In contrast, the analysis of hemivariational inequalities uses as main ingredient the properties 
of the subdifferential in the sense of Clarke \cite{Clarke}, defined for locally Lipschitz functions, 
which may be nonconvex. 
Hemivariational inequalities were first introduced in early eighty's by Panagiotopoulos in the context of 
applications in engineering problems.  Studies of variational and hemivariational inequalities can be 
found in several comprehensive references, e.g., \cite{BC, Br,CWW,Gl,GLT,LXH,NP,WXWC,XHC,XHW,ZHJ}. 
Variational-hemivariational inequalities represent 
a special class of inequalities, in which both convex and nonconvex functions are present. A recent reference in the field is the monograph \cite{SMBOOK}. There, existence, uniqueness and convergence results have been obtained and then used in the study of various mathematical models which describe the contact between a deformable body and a foundation.

The optimal control theory deals with the existence and, when possible, the uniqueness of optimal pairs and optimal control. It also deals with the derivation of necessary conditions of optimality or, better, necessary and   sufficient conditions of optimality. 
Optimal control problems for variational and hemivariational inequalities  have been discussed in several works, including \cite{Ba,F,LZ, MM2, MP, NST,SX1,SX2,Ti,Tiba}. Due to the nonsmooth and nonconvex feature of the functionals involved, the treatment of optimal control problems for such inequalities requires the use of their approximation by smooth optimization problems. And, on this matter, establishing convergence results for the optimal pairs represents a topic  of major interest.

In~\cite{S1} we have studied variational-hemivariational inequalities
of the form: find $u\in K_g$ such that
\begin{equation}\label{1}
\langle A u, v - u \rangle + \varphi (u, v) - \varphi (u, u) + j^0(u; v - u)
\ge (f,\pi v-\pi u)_Y \quad\forall\,v \in K_g.
\end{equation}

\noindent
Here and everywhere in this paper $X$ is a reflexive Banach space,   $\langle\cdot,\cdot\rangle$  denotes the duality pairing
between $X^*$ and its dual $X^*$, $Y$ is a real Hilbert space endowed with the inner product  $(\cdot,\cdot)_Y$ and $\pi:X\to Y$.  Moreover, in (\ref{1}) we supposed that $K_g=gK$  where $K\subset X$ and $g>0$,
$A \colon X \to X^*$,  $\varphi \colon X \times X \to \real$,
$j \colon X \to \real$ and $f\in Y$.
Note that the function
$\varphi(u, \cdot)$ is assumed to be convex
and the function $j$ is locally Lipschitz and, in general, nonconvex. Therefore, following the terminology in the literature, 
we see that inequality (\ref{1}) represents a
variational-hemivariational inequality.

A short description of the results obtained in \cite{S1} is the following. First,
the existence and uniqueness of the solution of (\ref{1}) was proved  by using a result proved in  \cite{MOS30}, based
on arguments of surjectivity for pseudomonotone operators
and the Banach fixed point argument. Then,  under specific assumptions on the functions $\varphi$ and $j$,  the continuous dependence
of the solution with respect $f$ and $g$ was studied  and a convergence result was  proved. Next,  two optimal control problems were considered, in which the control were $f$ and $g$, respectively. The existence of  optimal pairs together with some convergence results were proved, for each problem.
Finally, these abstract results were used in the study of a one-dimensional mathematical model which describes the equilibrium
of an elastic rod in unilateral contact with a foundation, under the action of a body force.

As it results from above, the study in \cite{S1} was focused on the dependence of the  solution with respect the parameters $f$ and $g$ or, equivalently, with respect to $f$ and the set of constraints $K_g$.
Nevertheless, various examples can be considered in which the solution depends on a number of parameters which appear in the operator $A$, or in the functions $\varphi$ and $j$, as well. All of these parameters could play the role of control in optimal control problems associated to inequality (\ref{1}).
For this reason, there is a need to extend the results in \cite{S1} to more general cases and this is the aim of this current paper. Considering such extension leads to various mathematical difficulties. To overcome them we use new assumptions and new arguments, different to those used in \cite{S1}, which represents the trait of novelty of this current paper. 

The rest of the manuscript is organized as follows. In Section~\ref{s2}
we review some notation  and preliminary results. In Section~\ref{s3} we introduce a variational-hemivaria\-tional inequality in which all the data depend on a parameter $p$. We state the behavior of the solution of this inequality with respect to $p$ and provide a convergence result. Then, in Section
\ref{s4} we consider a class of optimal control problems associated to the variational-hemivariational inequality, for which we prove
the existence and convergence of the optimal pairs. Finally,
in Section~\ref{s5} we give an example which
illustrate a potential application of our abstract study.  The example arises from Contact Mechanics and is given by a variational-hemivariational inequality which describes the contact of an elastic body with a foundation made of a rigid body covered by a layer of deformable material.

\section{Preliminaries}\label{s2}
\setcounter{equation}0

We use notation $\|\cdot\|_X$ and $0_X$ for the norm and the zero space element of $X$, respectively.
Unless stated otherwise, all the limits, upper and lower limits below are considered as $n\to\infty$, even if we do not mention it explicitly. The symbols ``$\rightharpoonup$"  and ``$\to$" 
denote the weak and the strong convergence in various spaces which will be specified. Nevertheless, for simplicity, we write $g_n\to g$ for the convergence in $\mathbb{R}$.

\begin{Definition}
An operator $A \colon X \to X^*$ is said to be:
	
	{\rm a)  monotone},
	if for all $u$, $v \in X$, we have $\langle Au - A v, u-v \rangle \ge 0$;

	{\rm b) bounded}, if $A$ maps bounded sets of $X$
	into bounded sets of $X^*$;

	{\rm c) pseudomonotone},
	if it is bounded and $u_n \to u$ weakly in $X$ with
	$$\displaystyle \limsup\,\langle A u_n, u_n -u \rangle \le 0$$
	implies\ \ $\displaystyle \liminf\, \langle A u_n, u_n - v \rangle\ge \langle A u, u - v \rangle$\ for all $v \in X$.
\end{Definition}

We now recall the definition of the Clarke subdifferential for a locally Lipschitz function.

\begin{Definition}
A function $j \colon X \to \real$ is said to be 
		{\rm locally Lipschitz}
		if for every 
		$x \in X$ there exists $U_x$ a neighborhood of $x$ and a constant $L_x>0$
		such that 
		\[
		|j(y) - j(z)| \le L_x \| y - z \|_X\qquad{\rm 	for\ all\ \ }\ y,\ z \in U_x.\]
	 The
	{\rm generalized (Clarke) directional derivative} of $j$ at the point
	$x \in X$ in the direction $v \in X$ is defined
	by
	\begin{equation*}
	j^{0}(x; v) = \limsup_{y \to x, \ \lambda \downarrow 0}
	\frac{j(y + \lambda v) - j(y)}{\lambda}.
	\end{equation*}
	The {\rm generalized gradient (subdifferential)} of $j$ at $x$
	is a subset of the dual space $X^*$ given by
	\begin{equation*}
	\partial j (x) = \{\, \zeta \in X^* \mid j^{0}(x; v) \ge
	{\langle \zeta, v \rangle} \quad \forall\, v \in X \, \}.
	\end{equation*}
	A locally Lipschitz function $j$ is said to be {\rm regular} (in the sense of Clarke) at the point $x\in X$ if for all $v\in X$ the one-sided directional derivative $j'(x;v)$ exists and $j^0(x;v)=j'(x;v)$.
	
\end{Definition}

We shall use the following properties of the generalized directional derivative and the generalized gradient.

\begin{Proposition}\label{subdiff}
	Assume that $j \colon X \to \real$ is a locally Lipschitz function. Then the following hold:
	
	\medskip{\rm (i)} 
	For every $x \in X$, the function
	$X \ni v \mapsto h^0(x;v) \in \real$ is positively  homogeneous and subadditive, 
	i.e., $j^0(x; \lambda v) = \lambda j^0(x; v)$ for all
	$\lambda \ge 0$, $v\in X$ and $j^0 (x; v_1 + v_2) \le
	j^0(x; v_1) + j^0(x; v_2)$ for all $v_1$, $v_2 \in
	X$, respectively.
	
	\medskip{\rm (ii)}
	For every $v \in X$, we have
	$j^0(x; v) = \max \, \{ \, \langle \xi, v \rangle
	\mid \xi \in \partial j(x) \, \}$.
\end{Proposition}

We now recall the notion of convergence in the sense of Mosco, denoted  by $``\stackrel{M}{\longrightarrow}"$, which will be used in the Sections \ref{s3}--\ref{s5} of this paper.

\begin{Definition}\label{M} Let $X$ be a normed space, $\{K_n\}$  a sequence of nonempty subsets of $X$ and $K$  a nonempty subset of $X$. We say that the sequence $\{K_n\}$ converges to $K$ in the Mosco sense if the following conditions hold.

	\bigskip\noindent
	$(M_1)$\qquad$\left\{
	\begin{array}{l}
	\mbox{For each $v\in K$ there exists a sequence $\{v_n\}$ such that}\\[1mm]
	\mbox{$v_n\in K_n$ for each $n\in \mathbb{N}$  and $v_n\to v$ in $X$.} 
	\end{array}\right.$
	
	\bigskip\noindent
	$(M_2)$\qquad$\left\{
	\begin{array}{l}
	\mbox{For each sequence $\{v_n\}$ such that}\\[1mm]
	\mbox{$v_n\in K_n$ for each $n\in \mathbb{N}$  and $v_n\rightharpoonup v$ in $X$ we have}\ v\in K. 
	\end{array}\right.$
	
\end{Definition}

Note that the convergence in the sense of Mosco depends on the topology of the normed space $X$ and, for this reason,  we write it  explicitly
$K_n\stackrel{M}{\longrightarrow} K_n$ in $X$. 
More details on this topic could be found in \cite{Mosco}.

\medskip
We proceed with the following version of the Weierstrass theorem.

\begin{Theorem}\label{tW}
	Let $X$ be a reflexive Banach space, $K$  a nonempty weakly closed 
	subset of $X$ and   $J:X\to\mathbb{R}$ a weakly lower semicontinuous function. In addition, assume that either $K$ is bounded or $J$ is coercive,  
	i.e., $J(v)\to\infty$ as $\|v\|_X\to\infty$.
	Then, there exists at least an element $u$ such that 
	\begin{equation}\label{04}
	u\in K,\qquad J(u)\le J(v)\qquad\forall\, v\in K.
	\end{equation}
\end{Theorem}

The proof of Theorem \ref{tW}  can be found in many books and survey, see, for instance, \cite{Ku,SofMat}.

\medskip

We now  introduce the variational-hemivariational inequality we are interested in and state its unique solvability.
The functional framework is the following. Let $Z$ be a normed space and let $\Lambda\subset Z$ be a set of parameters.  For any 
parameter $p\in\Lambda$ we consider
a set $K_p\subset X$, an operator
$A_p \colon X \to X^*$,  the functions $\varphi_p \colon X \times X \to \R$ and
$j_p\colon X \to \R$ and an element $f_p\in Y$. With these data, for $p\in\Lambda$ given, we consider the following inequality problem.

\medskip\noindent{\bf Problem}  ${\cal P}$. {\it  Find an element  $u=u(p)$ such that}
\begin{eqnarray}
&&\label{2}u\in K_p,\quad \langle A_pu, v - u \rangle + \varphi_p (u,v)-\varphi_p(u,u) + j^0_p(u; v - u)\\[2mm]
&&\qquad\qquad\qquad\ge (f_p, \pi v - \pi u)_Y \qquad\forall\,  v \in K_p.\nonumber
\end{eqnarray}


\medskip
In the study of Problem ${\cal P}$, we consider the following
hypotheses on the data.

\begin{eqnarray}\label{K}
&&\quad\ K_p \ \mbox{is nonempty, closed and convex subset of} \ X.\\[6mm]
&&\label{A}
\left\{
\begin{array}{l}
A_p \colon X \to X^* \ \mbox{is such that} \\[2mm]
\ \ {\rm (a)} \ \mbox{it is pseudomonotone}; 
\\ [2mm]
\ \ {\rm (b)} \ \mbox{it is strongly monotone, i.e., there exists}\ m_p > 0  \ \mbox{such that}\\[2mm]
\qquad\quad \langle Av_1 - Av_2, v_1 - v_2 \rangle
\ge
m_p \| v_1 - v_2 \|_X^{2} \quad \mbox{for all} \ v_1, v_2 \in X.
\end{array}\right.
\\[6mm]
&&\left\{\begin{array}{l}
\varphi_p \colon X \times X \to \real \ \mbox{is such that}\\ [2mm]
\ \ {\rm (a)} \
\varphi_p(u, \cdot) \colon X \to \real \ \mbox{is convex and lower semicontinuous},\\[2mm]
\qquad \ \ \mbox{for all} \ u \in X; \\ [2mm]
\ \ {\rm (b)} \
\mbox{there exists} \ \alpha_p > 0 \ \mbox{such that}  \\ [2mm]
\qquad 
\varphi_p(u_1, v_2) - \varphi_p(u_1, v_1) + \varphi_p(u_2, v_1) - \varphi_p(u_2, v_2)\\ [2mm]
\qquad\quad\le \alpha_p \| u_1 - u_2 \|_X \, \| v_1 - v_2 \|_X\ \ 
\ \mbox{for all} \ u_1, u_2, v_1, v_2 \in X.
\end{array}
\right.
\label{fi}
\end{eqnarray}
\begin{equation}
\label{j} \left\{
\begin{array}{l}
j_p \colon X \to \real \ \mbox{is such that}\\ [2mm]
\ \ {\rm (a)} \
\mbox{it is locally Lipschitz;} \\ [2mm]
\ \ {\rm (b)} \
\| \xi \|_{X^*} \le c_{0p} + c_{1p} \, \| v \|_X \ \mbox{for all}
\ v \in X,\ \xi\in \partial j_p(v),\qquad \\ [2mm]
\qquad\quad \mbox{with} \ c_{0p}, c_{1p} \ge 0; \\ [2mm]
\ \ {\rm (c)} \
\mbox{there exists} \ \beta_p > 0 \ \mbox{such that} \\ [2mm]
\qquad \quad
j^0_p(v_1; v_2 - v_1) + j^0_p(v_2; v_1 - v_2) \le \beta_p \, \| v_1 - v_2 \|_X^2 \\ [2mm]
\qquad \quad \mbox{for all} \ v_1, v_2 \in X.
\end{array}
\right.
\end{equation}
\begin{eqnarray}
&&\label{sm}
\qquad\quad\alpha_p + \beta_p < m_p.\\[2mm]
&&\label{f}\qquad\quad f_p\in Y.\\ [4mm]
&&\label{pi}
\left\{ \begin{array}{l} \pi \colon X \to Y\ {\rm is\ a\ linear\ continuous\ operator,\ i.e.,}\qquad\qquad\\[2mm]
\mbox\ \|\pi v\|_{Y} \leq d_0\,\|v\|_X\quad \forall\,v\in X\ \ {\rm with}\ d_0>0.\\[2mm]
\end{array}\right. 
\end{eqnarray}

\medskip

The unique solvability of the variational-hemivariational inequality (\ref{1}) is provided by the following result.

\begin{Theorem}\label{t0}
	Assume  $(\ref{K})$--$(\ref{pi})$.
	Then, Problem~ ${\cal P}$ has a unique solution $u\in K_p$.
\end{Theorem}

Theorem~\ref{t0} represent a slightly modified version of Theorem 18 in \cite{MOS30}. Its proof is
is carried out in several steps, by using the properties of the subdifferential, both in the sense of Clarke and in the sense of convex analysis, surjectivity result for  pseudomonotone multivalued operators, and the Banach fixed point argument. For additional details  we also refer the reader to \cite{S1} and \cite[Remark 13]{SMBOOK}.

\section{A convergence result}\label{s3}
\setcounter{equation}0

\medskip Theorem \ref{t0} allows us to define the operator $p\mapsto u(p)$ which associates to each  $p\in\Lambda$ the solution $u=u(p)\in K_p$ of the variational-hemivaria\-tional inequality (\ref{2}). An important property of this operator is its  continuity in a specific sense, which we state and prove in this section, under additional assumptions. This property represents a crucial ingredient in the study of optimal control problems associated to inequality (\ref{2}). 

Let $\{p_n\}$ be a sequence of elements in $\Lambda$ and assume that, for each $n\in\mathbb{N}$, the corresponding set $K_{p_n}$, operator $A_{p_n}$, functions $\varphi_{p_n}$, $j_{p_n}$ and element $f_{p_n}$ satisfy assumptions (\ref{K})--(\ref{f}) with constants $m_{p_n}$, $\alpha_{p_n}$, $c_{0p_n}$ $c_{1p_n}$, $\beta_{p_n}$. To avoid any confusion, when used for $p=p_n$, we refer to these assumptions as assumptions (\ref{K})$_n$--(\ref{f})$_n$.
Then, if condition (\ref{pi}) is satisfied, we deduce from Theorem \ref{t0} that for each $n\in\mathbb{N}$ there exists a unique solution $u_n=u(p_n)$ for the following problem.

\medskip\noindent{\bf Problem}  ${\cal P}_n$. {\it Find an element  $u_n=u(p_n)$ such that}
\begin{eqnarray}
&&\hspace{-10mm}\label{2n}u_n\in K_{p_n},\quad \langle A_{p_n}u_n, v - u_n \rangle + \varphi_{p_n} (u_n,v)-\varphi_{p_n}(u_n,u_n) + j^0_{p_n}(u_n; v - u_n)\\[2mm]
&&\hspace{-10mm}\qquad\qquad\qquad\ge (f_{p_n}, \pi v - \pi u_n)_Y \quad\forall\,  v \in K_{p_n}.\nonumber
\end{eqnarray}

\medskip
We now  consider the following additional assumptions.
\begin{eqnarray}
&&\label{cvK}\quad\ K_{p_n}\stackrel{M}{\longrightarrow} K_p\quad {\rm in}\quad X,\quad{\rm as}\quad n\to\infty.
\\ [3mm]
&&\label{cvA}
\left\{
\begin{array}{ll}
\mbox{For any}\ n\in\mathbb{N}\ \mbox{there exists }F_n\ge 0\ \mbox{and}\ \delta_n\ge 0\  \mbox{ such
	that  
}\\[2mm]
\mbox{(a) } \|A_{p_n}v-A_{p}v\|_X\le F_n(\|v\|_X+\delta_n)\quad\mbox{for all}\ v\in X;
\\[3mm]
\mbox{(b) }\displaystyle\lim_{n\to0}F_n=0;\\[3mm]
\mbox{(c) the sequence}\ \{\delta_n\}\subset\mathbb{R}\ \mbox {is bounded}.
\end{array}
\right.\\ [3mm]
&&\label{cvfi}\left\{\begin{array}{ll} 
\mbox{For any}\ n\in\mathbb{N}\  \mbox{there exists a function}\ c_n:\real_+\to \real_+\ \mbox{such that}
\\ [2mm]
\mbox{(a) }
\varphi_{p_n}(u,v_1)-\varphi_{p_n}(u,v_2)\le c_n(\|u\|_X){\|v_1-v_2\|_X}\\ [2mm]
\qquad\forall\,u,\,v_1,\,v_2\in
X;
\\ [2mm]
\mbox{(b)  the sequence}\ \{c_n(\|v_n\|_X)\}\subset\mathbb{R}\ \mbox {is bounded}\\ [2mm]
\qquad \mbox{whenever the sequence} \ \{v_n\}\subset X\ \mbox {is bounded}.
\end{array}\right.
\\ [3mm]
&&\label{cvfin}
\left\{\begin{array}{ll}\mbox{For all sequences}\ \{u_n\},\ \{v_n\}\subset X\  \mbox{such that}\\[2mm]
u_n\rightharpoonup u\ \mbox{in}\ X,\ v_n\to v\ \mbox{in}\ X,\ \mbox{we have}\hspace{26mm} \\ [2mm]
\limsup\,\big(\varphi_{p_n}(u_n,v_n)-\varphi_{p_n}(u_n,u_n)\big)\le \varphi_p(u,v)-\varphi_p(u,u).
\end{array}\right.\\ [3mm]
&&\label{cvj}
\left\{\begin{array}{ll}\mbox{For all sequences}\ \{u_n\},\ \{v_n\}\subset X\  \mbox{such that}\\[2mm]
u_n\rightharpoonup u\ \mbox{in}\ X,\ v_n\to v\ \mbox{in}\ X,\  \mbox{we have} \\ [2mm]
\ \limsup j^0_{p_n}(u_n; v_n - u_n) \le j^0_p(u; v - u).
\end{array}
\right.\\[3mm]
&&\label{cvpi}
\left\{\begin{array}{ll}\mbox{For all sequence}\ \{v_n\}\subset X\  \mbox{such that}\\[2mm]
v_n\rightharpoonup v\ \mbox{in}\ X,\  \mbox{we have}
\quad \pi v_n\to\pi v\quad{\rm in}\quad Y.
\end{array}\right.\\[3mm]
&&\label{cvf}\quad\ f_{p_n}\rightharpoonup f_p\quad {\rm in}\quad Y,\quad{\rm as}\quad n\to\infty.
\\[3mm]
&&\label{cvs}\mbox{There exists}\ m_0>0\ \mbox{such that}\ m_0+\alpha_{p_n}+\beta_{p_n}\le m_{p_n}\ \ \forall\, n\in\mathbb{N}.
\\[3mm]
&&\label{cvb}\mbox{The sequences}\ \{c_{0p_n}\},\, \{c_{1p_n}\}\subset \real\ \mbox{are bounded}.
\end{eqnarray}


\medskip
The main result of this section is the following.

\begin{Theorem}\label{t1}  Assume  $(\ref{K})$--$(\ref{pi})$ and, for each $n\in\mathbb{N}$, assume $(\ref{K})_n$--$(\ref{f})_n$. Moreover, assume $(\ref{cvK})$--$(\ref{cvb})$ and denote by $u_n$ and $u$ the solutions of Problems ${\cal P}_n$ and ${\cal P}$, respectively. Then
the sequence $\{u_n\}$ converges strongly to $u$, i.e., 
\begin{equation}\label{cvu}
u_n\to u\qquad{\rm in}\quad X.
\end{equation}
\end{Theorem}

 Note  that in applications, conditions (\ref{cvK})--(\ref{cvj}) and (\ref{cvf})--(\ref{cvb})
are satisfied when $p_n\rightharpoonup p$ in $Z$. In this case Theorem \ref{t1} implies that
\begin{equation}\label{es}
p_n\rightharpoonup p\quad{\rm in} \quad Z\ \Longrightarrow\ u(p_n)\to u(p)\quad{\rm in}\quad X,
\end{equation}
which shows that the operator $p\mapsto u(p):\Lambda\to X$ is 
weakly-strongly continuous.

\medskip

The proof of Theorem \ref{t1} will be carried out
in several steps that we present in what follows.  Everywhere below we assume that the hypotheses of Theorem \ref{t1} hold. The first step of the proof is the following.

\begin{Lemma}\label{l1}	There is an element ${\widetilde{u}} \in K_p$ and
a subsequence of $\{ u_n \}$, 
still denoted by $\{ u_n \}$,
such that $u_n \to {\widetilde{u}}$ weakly in $X$, as $n\to \infty$.
\end{Lemma}	

\noindent{\bf Proof.}  We first establish the boundedness of $\{ u_n \}$ in $X$. Let $v$ be a given element in $K_p$. We use assumption (\ref{cvK})  and consider a sequence
$ \{v_n\}$ such that $v_n\in K_{p_n}$ for all $n\in\mathbb{N}$ and
\begin{equation}\label{vv}
v_n\to v\quad \mbox{in}\quad X.
\end{equation}
\noindent

Let $n\in\mathbb{N}$. We test with $v_n\in K_{p_n}$ in (\ref{2n}) to obtain
\begin{eqnarray}
&&\label{5}
\langle A_{p_n}u_n, u_n - v_n \rangle\le
  \\ [2mm] 
&& 
\quad +
\varphi_{p_n}(u_n,v_n) - \varphi_{p_n}(u_n,u_n)+ \, j^0_{p_n}(u_n; v_n - u_n) +
(f_{p_n}, \pi u_n - \pi v_n)_Y. \nonumber
\end{eqnarray}
Moreover, using the strongly monotonicity of the operator $A_{p_n}$ we  deduce that
\begin{equation}\label{6}
\langle A_{p_n}u_n, u_n - v_n \rangle\ge m_{p_n}\| u_n  - v_n \|_X^2 +\langle A_{p_n}v_n, u_n - v_n \rangle.
\end{equation}
Next, we write
\begin{eqnarray*}
&&
\varphi_{p_n}(u_n,v_n) - \varphi_{p_n}(u_n,u_n)\\ [2mm]
&&\quad	=\varphi_{p_n}(u_n,v_n) - \varphi_{p_n}(u_n,u_n)+ 	\varphi_{p_n}(v_n,u_n) - \varphi_{p_n}(v_n,v_n)\nonumber\\ [2mm]
&&\qquad+\varphi_{p_n}(v_n,v_n) - \varphi_{p_n}(v_n,u_n),\nonumber
\end{eqnarray*}
then we use assumptions (\ref{fi})$_n$(b), (\ref{cvfi}) to see that
\begin{eqnarray}
&&\label{7}
\varphi_{p_n}(u_n,v_n) - \varphi_{p_n}(u_n,u_n)\\ [2mm]
&&\qquad\le \alpha_{p_n}\| u_n  - v_n \|_X^2+c_n(\| v_n\|_X) \| u_n  - v_n \|_X.\nonumber
\end{eqnarray}
On the other hand, by (\ref{j})$_n$(b) and Proposition~\ref{subdiff}(ii), we have
\begin{eqnarray*}
&& j^0_{p_n}(u_n; v_n - u_n) \\[3mm]
&&\quad = j^0_{p_n}(u_n; v_n - u_n) + j^0_{p_n} (v_n; u_n - v_n) - j^0_{p_n} (v_n; u_n - v_n)\nonumber \\ [3mm]
&&
\qquad \le j^0_{p_n}(u_n; v_n - u_n) + j^0_{p_n} (v_n; u_n - v_n) 
+ | j^0_{p_n} (v_n; u_n - v_n)| \nonumber \\ [3mm]
&&\quad \qquad
\le \beta_{p_n} \| u_n - v_n \|_X^2 +
|\max \{ \langle \zeta, u_n - v_n \rangle \mid \zeta \in \partial j_{p_n}(v_n) \} | \nonumber 
\end{eqnarray*}
and, using assumption (\ref{j})$_n$(b) we find that
\begin{eqnarray}\label{8}
&& j^0_{p_n}(u_n; v_n - u_n)\le\beta_{p_n} \| u_n - v_n \|_X^2 + (c_{0p_n} + c_{1p_n} \| v_n \|_X) \| u_n - v_n \|_X. 
\end{eqnarray}
Finally, assumption (\ref{pi}) yields
\begin{equation}\label{9}
(f_{p_n}, \pi u_n - \pi v_n)_Y  \le d_0\| f_{p_n} \|_Y\| u_n - v_n \|_X.
\end{equation}

We now combine inequalities (\ref{5})--(\ref{9}) to see that
\begin{eqnarray*}
&&
m_{p_n} \, \| u_n  - v_n \|_X^2 \le \langle A_{p_n}v_n, v_n - u_n \rangle\\ [2mm]
&&\quad+\alpha_{p_n}\| u_n  - v_n \|_X^2+c_n(\| v_n\|_X) \| u_n  - v_n \|_X\nonumber
\\ [2mm]
&& 
\qquad +\beta_{p_n} \| u_n - v_n \|_X^2 + (c_{0p_n} + c_{1p_n} \| v_n \|_X) \| u_n - v_n \|_X+d_0\| f_{p_n} \|_Y\| u_n - v_n \|_X,\nonumber
\end{eqnarray*}
which implies that
\begin{eqnarray}
	&&\label{10}
	(m_{p_n}-\alpha_{p_n}-\beta_{p_n}) \, \| u_n  - v_n \|_X \le \|A_{p_n}v_n\|_X\\ [2mm]
	&&\quad+c_n(\| v_n\|_X)+(c_{0p_n} + c_{1p_n} \| v_n \|_X)+d_0\| f_{p_n} \|_Y.\nonumber
\end{eqnarray}
Note that assumptions (\ref{cvs}) and (\ref{cvA}) imply that 
\[m_0\le m_{p_n}-\alpha_{p_n}-\beta_{p_n}\quad\mbox{and}\quad\|A_{p_n}v_n\|_X\le \|A_{p}v_n\|_X+F_n(\|v_n\|_X+\delta_n),\]
respectively. Therefore, using these inequalities in (\ref{10}) we deduce that
\begin{eqnarray}
&&\label{10n}
m_0 \, \| u_n  - v_n \|_X\le \|A_{p}v_n\|_X+F_n(\|v_n\|_X+\delta_n)\\ [2mm]
&&\quad+c_n(\| v_n\|_X)+(c_{0p_n} + c_{1p_n} \| v_n \|_X)+d_0\| f_{p_n} \|_Y.\nonumber
\end{eqnarray}

Next, by (\ref{vv}) and (\ref{cvf}) we know that the sequences $\{v_n\}$ and $\{f_{p_n}\}$ are bounded in $X$ and $Y$, respectively. Therefore, using
assumptions (\ref{A})(a), (\ref{cvA})(b), (c),  (\ref{cvfi})(b),  (\ref{cvb}) and  inequality (\ref{10n}),
we deduce that  that there is a constant $C > 0$
independent of $n$ such that $\| u_n-v_n\|_X \le C$. 
This implies that  $\{u_n\}$ is a bounded sequence in $X$ and,
from the reflexivity of $X$, by passing to
a subsequence, if necessary,  we deduce that
\begin{equation}\label{10m}
u_n \rightharpoonup {\widetilde{u}} \ \ \mbox{in} \ \ X, \ \mbox{as} \ n\to \infty
\end{equation}
with some ${\widetilde{u}} \in X$.  Finally, recall that $u_n\in K_{p_n}$. Therefore, the convergence (\ref{10m}) combined with assumption (\ref{cvK}) shows that
$\widetilde{u}\in K_p$ which concludes the proof. \hfill$\Box$

\medskip
The second step in the proof is as follows.

\begin{Lemma}\label{l2}	The element ${\widetilde{u}} \in K_p$ is a solution of
	Problem ${\cal P}$.
\end{Lemma}	

\noindent{\bf Proof.}	Let $v$ be a given element in $K_p$. We use assumption (\ref{cvK})  and consider a sequence
$ \{v_n\}$ such that $v_n\in K_{p_n}$ for all $n\in\mathbb{N}$ and (\ref{vv}) holds.
Let $n\in\mathbb{N}$. Then inequality (\ref{5}) holds and, passing to the upper limit in this inequality, we find that
\begin{eqnarray}
&&\label{20}
\hspace{-5mm}\limsup\,\langle A_{p_n}u_n, u_n - v_n \rangle\le 
\limsup\,[\varphi_{p_n}(u_n,v_n) - \varphi_{p_n}(u_n,u_n)]
\\ [2mm] 
&& 
\hspace{-2mm}\quad + \limsup\,j^0_{p_n}(u_n; v_n - u_n) +
\limsup\, (f_{p_n}, \pi u_n - \pi v_n)_Y. \nonumber
\end{eqnarray}
We now use the convergences (\ref{10m}), (\ref{vv}) and assumptions (\ref{cvfin}), (\ref{cvj}), to deduce that
\begin{eqnarray}
&&\label{21}\limsup\,[\varphi_{p_n}(u_n,v_n) - \varphi_{p_n}(u_n,u_n)]\le \varphi_{p}(\widetilde{u},v) - \varphi_{p}(\widetilde{u},\widetilde{u}),\\ [2mm]
&&\label{22}\limsup\,j^0_{p_n}(u_n; v_n - u_n) \le j^0_{p}(\widetilde{u};v-\widetilde{u}). 
\end{eqnarray}
On the other hand, the convergences (\ref{10m}), (\ref{vv}) and assumptions (\ref{cvpi}), (\ref{cvf}) yield
\begin{equation}\label{23}
\lim\, (f_{p_n}, \pi u_n - \pi v_n)_Y = (f_{p}, \pi \widetilde{u}- \pi v)_Y.
\end{equation}
We now combine relations (\ref{20})--(\ref{23}) to deduce that
\begin{eqnarray}
&&\label{24}
\limsup\,\langle A_{p_n}u_n, u_n - v_n \rangle\le 
\varphi_{p}(\widetilde{u},v) - \varphi_{p}(\widetilde{u},\widetilde{u})
\\ [2mm] 
&& \qquad\qquad + j^0_{p}(\widetilde{u};v-\widetilde{u})+ (f_{p}, \pi \widetilde{u}- \pi v)_Y. \nonumber
\end{eqnarray}

Next, we write
\begin{eqnarray*}
&&
\langle A_{p_n}u_n, u_n - v_n \rangle= \langle A_{p}u_n, u_n - v \rangle+\langle A_{p}u_n, v - v_n \rangle+\langle A_{p_n}u_n-A_pu_n, u_n - v_n \rangle.
\end{eqnarray*}
Then we use we assumptions (\ref{A})(a), (\ref{cvA}) and the convergence (\ref{vv}) to see that
\begin{eqnarray*}
&&\langle A_{p}u_n, v - v_n \rangle\to 0, \\ [2mm]
&&
|\langle A_{p_n}u_n-A_pu_n, u_n - v_n \rangle|\le F_n(\|u_n\|_X+\delta_n)\|u_n-v_n\|_X\to 0,
\end{eqnarray*}
which imply that
\begin{equation}\label{25}
\limsup\,\langle A_{p_n}u_n, u_n - v_n \rangle= \limsup\, \langle A_{p}u_n, u_n - v\rangle. 
\end{equation}

We now combine inequality (\ref{24}) with  (\ref{25}) to obtain that
\begin{eqnarray}
&&\label{26}
\limsup\,\langle A_{p}u_n, u_n - v\rangle\le 
\varphi_{p}(\widetilde{u},v) - \varphi_{p}(\widetilde{u},\widetilde{u})
\\ [2mm] 
&& \qquad\qquad + j^0_{p}(\widetilde{u};v-\widetilde{u})+ (f_{p}, \pi \widetilde{u}- \pi v)_Y. \nonumber
\end{eqnarray}
Next, we take $v=\widetilde{u}$ in (\ref{26}) and use Proposition \ref{subdiff}(i) to deduce that
\begin{eqnarray}
&&\label{26n}
\limsup\,\langle A_{p}u_n, u_n -\widetilde{u} \rangle\le 0.
\end{eqnarray}
Exploiting now the pseudomonotonicity of $A_p$, 
from (\ref{10m}) and (\ref{26n}), we have
\begin{equation}\label{27}
\langle A_p{\widetilde{u}}, {\widetilde{u}} - v \rangle \le
\liminf\, \langle A_p u_n, u_n - v \rangle \ \ \mbox{for all} \ \ v \in X.
\end{equation}
Next, from  (\ref{27}) and (\ref{26}) we obtain that $\widetilde{u}$ satisfies the inequality (\ref{2}), which concludes the proof.
\hfill$\Box$

\medskip We are now in position to provide the proof of Theorem \ref{t1}.

\medskip \noindent{\bf Proof.}	Since Problem~${\cal P}$ has a unique solution, denoted $u$, it follows from Lemma \ref{l2} that ${\widetilde{u}} = u$. This implies that every subsequence of $\{ u_n \}$
which converges weakly in $X$ has the same limit and, therefore, it follows that the whole sequence $\{ u_n \}$ converges weakly in $X$ to $u$,  as $n\to \infty$.

We now prove that $u_n \to u$ in $X$, as $n\to\infty$. To this end, we take 
	$v = {\widetilde{u}} \in K_p$ in both (\ref{27}) and (\ref{26}), then we use  equality ${\widetilde{u}} = u$ to obtain
	\[
	0 \le \liminf\, \langle A_pu_n, u_n - u \rangle
	\le
	\limsup\,\langle A_pu_n, u_n - u \rangle \le 0,
	\]
	which shows that $\langle A_pu_n, u_n - u \rangle \to 0$,
	as $n\to\infty$.
	Therefore, using the strong monotonicity of the operator $A_p$
	and the convergence $u_n\rightharpoonup u$ in $X$,
	we have
	\[
	m_A \| u_n - u \|_X^2 \le
	\langle A_pu_n - A_pu, u_n - u \rangle =
	\langle A_pu_n, u_n - u \rangle - \langle A_pu, u_n - u \rangle \to 0,
	\]
	as $n\to \infty$. Hence, it follows that $u_n \to u$ in $X$, which concludes the proof of the theorem.
\hfill$\Box$

\section{An optimal control problem}\label{s4}
\setcounter{equation}0

Below in this section, we assume that the  normed space $Z$ has a particular structure, i.e., it is of the form
$Z =Q\times \Theta$, 
where $Q$ is a reflexive Banach space with the norm $\|\cdot\|_{Q}$ and $\Theta$ is a normed space endowed with the norm $\|\cdot\|_{\Theta}$.  A generic element of $Z$ will be denoted by $p=(q,\eta)$. 
We endow $Z$ with the norm
\[\|p\|_{Z}=\|q\|_{Q}+\|\eta\|_{\Theta}\qquad\forall\, p=(q,\eta)\in Z.\]

Let $U\subset Q$ and $\Sigma\subset\Theta$ be given nonempty subsets assumed to be weakly closed in $Q$ and $\Theta$, respectively, and let $\Lambda=U\times \Sigma$. 
Then, for each $q\in U$ and
$\eta\in\Sigma$ we have $p=(q,\eta)\in\Lambda$ and, under the assumption of Theorem \ref{t0}, we denote by
$u=u(p)=u(q,\eta)$ the solution of Problem  ${\cal P}$.
Moreover, for each $\eta\in\Sigma$, let $F(\eta)$ denote a subset of $U$ which depends on $\eta$.  With these notation define
the set of admissible pairs for  Problem ${\cal P}$
by equality
\begin{equation}\label{vad}
{\cal V}_{ad}(\eta) = \{\,(u, q) \ : \  q\in F(\eta),\ u=u(q,\eta)\,\}.
\end{equation}
In other words, a pair $(u,q)$ belongs to
${\cal V}_{ad}(\eta)$ if and only if $q\in F(\eta)$ and, moreover, $u$ is the solution of  Problem ${\cal P}$ with $p=(q,\eta)$.
Consider also
a cost functional ${\cal L}:X\times U\to\mathbb{R}$. 
Then, the optimal control problem we are interested in is the following.

\medskip\noindent
{\bf Problem} ${\cal Q}$. {\it Given $\eta\in\Sigma$, find $(u^*, q^*)\in {\cal V}_{ad}(\eta)$ such that}
\begin{equation}\label{o1}
{\cal L}(u^*,q^*)=\min_{(u,q)\in {\cal V}_{ad}(\eta)} {\cal L}(u,q).
\end{equation}

\medskip
To solve Problem ${\cal Q}$ we consider the following assumptions.

\begin{eqnarray}
&&\label{o2} F(\eta)\  \ \mbox{is a nonempty weakly closed subset of}\ \  Q.\qquad
\\[4mm]
&&\label{o3}\left\{\begin{array}{l}
\mbox{For all sequences }\{u_k\}\subset X\mbox{ and
}\{q_k\}\subset
U \mbox{ such that}\\[2mm]
u_k\rightarrow u\mbox{\ \ in\ \ }X,\ \ q_{k}\rightharpoonup
q \mbox{\ \ in\ \ }
Q,\ \mbox{we have} \\[3mm]
{\rm (a)}\ \ \displaystyle\liminf_{k\to
	\infty}\,{\cal L}(u_k,q_k)\ge {\cal L}(u,q),
\\[3mm]
{\rm (b)}\ \ \displaystyle\lim_{k\to
	\infty}\,\big[{\cal L}(u_k,q_k)-{\cal L}(u,q_k)\big]=0.
\end{array}\right.
\\[4mm]
&&\label{o4}\left\{ \begin{array}{l}
\mbox{ There exists}\ h: U\to\R\ \mbox{such that}
\\[2mm]
{\rm (a)}\ \ {\cal L}(u,q)\ge h(q)\quad \forall\,u\in X,\  q\in U,\\ [2mm]
{\rm (b)}\ \ \|q_{k}\|_{Q}\to+\infty\ \Longrightarrow\ h(q_k)\to \infty.
\end{array}\right.\\[4mm]
&&\label{o5}\qquad U \ \ \mbox {is a bounded subset of} \ \  Q.
\end{eqnarray}

\begin{Example}\label{e0}
	A typical example of function ${\cal L}$ which satisfies conditions $(\ref{o3})$ and $(\ref{o4})$ is obtained by taking
	\[
	{\cal L}(u,p)=g(u)+h(q)\qquad\forall\, u\in X,\ q\in U,
	\]
	where $g:X\to\R$ is a continuous positive function and $h: U\to \R$ is a lower semicontinuous coercive function, i.e., it satisfies condition $(\ref{o4}){\rm (b)}$.		
\end{Example}

Our first result in this section is the following existence result.

\begin{Theorem}\label{t4}  Assume  $(\ref{K})$--$(\ref{f})$, for any $p\in\Lambda$. Moreover, assume $(\ref{pi})$, $(\ref{es})$, $(\ref{o2})$, $(\ref{o3})$  and, in addition, assume that either 
$(\ref{o4})$ or $(\ref{o5})$ hold. Then
Problem ${\cal Q}$ has  at least one solution $(u^*,q^*)$.
\end{Theorem}

The proof  will be carried out in several steps that we present in what follows. To this end, everywhere below we assume that the hypotheses of
Theorem \ref{t4} hold and, given $\eta\in\Sigma$, we  consider the function $J(\cdot,\eta):U\to\R$ defined by
	\begin{equation}\label{Jmn}
	J(q,\eta)=	{\cal L}(u(q,\eta),q)\qquad\forall\,q\in U.
	\end{equation}

The first step of the proof is as follows.

\begin{Lemma}\label{la} {\it For all sequences $\{q_k\}\subset U$ and
		$\{\eta_k\}\subset
		\Sigma$ such that
		$q_k\rightharpoonup q\mbox{ in }Q$,\,\,$\eta_k\rightharpoonup
		\eta \mbox{ in }
		\Theta$ { and for all  $s\in U$,
			the inequality below holds: }}
	\begin{equation}\label{j1}\displaystyle\limsup_{k\to
		\infty}\,[J(s,\eta_k)-J(q_k,\eta_k)]\leq J(s,\eta)-J(q,\eta).
	\end{equation}

	\end{Lemma}

\medskip\noindent{\bf Proof.}  Assume that $\{q_k\}\subset U$ and
$\{\eta_k\}\subset
\Lambda$ are two sequences  such that $q_k\rightharpoonup q\mbox{ in }Q$,\,\,$\eta_k\rightharpoonup
\eta\mbox{ in } \Theta$ and let $s\in U$. Then, using (\ref{Jmn}) we have
\begin{equation}\label{r1}
J(s,\eta_k)-J(q_k,\eta_k)={\cal L}(u(s,\eta_k),s)-{\cal L}(u(q_k,\eta_k),q_k).
\end{equation}
Moreover, assumption (\ref{es}) implies that   $u(s,\eta_k)\to u(s,\eta)$, $u(q_k,\eta_k)\to u(q,\eta)$, both in $X$. Therefore,
assumption (\ref{o3}) and definition (\ref{Jmn}) imply that
\begin{eqnarray}
&&\label{r2}
\lim_{k\to\infty} {\cal L}(u(s,\eta_k),s)={\cal L}(u(s,\eta),s)=J(s,\eta),\\[2mm]
&&\label{r3}\limsup_{k\to\infty}[-{\cal L}(u(q_k,\eta_k),q_k)]\le-{\cal L}(u(q,\eta),q)=-J(q,\eta).
\end{eqnarray}
We now pass to the upper limit in (\ref{r1}) and use relations (\ref{r2}), (\ref{r3}) to deduce that (\ref{j1}) holds.
\hfill$\Box$

\medskip
We now consider the following auxiliary problem.

\medskip\noindent
{\bf Problem} ${\cal R}$. {\it Given $\eta\in\Sigma$, find $q^*\in F(\eta)$ such that}
\begin{equation}\label{o1p}
J(q^*,\eta)=\min_{q\in F(\eta)} J(q,\eta).
\end{equation}

The solvability of Problem ${\cal R}$  is provided by the following result.

\begin{Lemma}\label{lb} 
Problem ${\cal R}$ has at least one solution $q^*$.
\end{Lemma}

\noindent{\bf Proof.}
Let $\eta\in\Sigma$. We take $\eta_k=\eta$ in the statement of Lemma \ref{la} to see that
for all sequences $\{q_k\}\subset U$ such that
$q_k\rightharpoonup q$ and for all $s\in U$
we have
\[\displaystyle\limsup_{k\to
	\infty}\,[J(s,\eta)-J(q_k,\eta)]\leq J(s,\eta)-J(q,\eta),\]
which implies that
\[\displaystyle\liminf_{k\to
	\infty}\,J(q_k,\eta)\ge J(q,\eta).\]
We conclude from here that the function $J(\cdot,\eta):U\to\mathbb{R}$ is lower semicontinuous.

Assume now that  (\ref{o4}) holds. Then,
for any sequence $\{q_k\}\subset U$,
we have
\[J(q_k,\eta)={\cal L}(u(q_k,\eta),q_k)\ge h(q_k).\]
Therefore, if $\|q_k\|_{Q} \to \infty$ we deduce that
$J(q_k,\eta)\to \infty$ which shows that the function $J(\cdot,\eta)$ is coercive.
Recall also the assumption (\ref{o2}) and the reflexivity of the space $Q$. The existence of at least one solution to Problem  ${\cal R}$ is now a direct consequence of Theorem \ref{tW}. 
On the other hand, if (\ref{o5}) holds we are still in position to apply  Theorem \ref{tW}. We deduce from here that, if either (\ref{o4}) or (\ref{o5})
hold, then Problem ${\cal R}$ has at least a solution, which concludes the proof.
\hfill$\Box$

\medskip
We are now in position to provide the proof of Theorem \ref{t4}.

\medskip\noindent{\bf Proof.} Using the definitions (\ref{Jmn}) and (\ref{vad}) it is easy to see that
\begin{equation}\label{b1}
\left\{ \begin{array}{l}
(u^*,q^*) \ \mbox{is a solution of Problem}\ {\cal Q}\ \mbox{if and only if}  \\ [2mm]
q^* \ \mbox{is a solution of Problem}\ {\cal R}\ \mbox{and}\ u^*=u(q^*,\eta).
\end{array}\right.  
\end{equation}	
Theorem \ref{t4} is a direct consequence of the equivalence (\ref{b1}) and Lemma \ref{lb}.
\hfill$\Box$

\medskip

\medskip

Let $\eta\in\Sigma$ and, for  each $n\in\mathbb{N}$,  consider a perturbation $\eta_n\in\Sigma$ of $\eta$ together with
the set of admissible pairs defined by
\begin{equation}\label{vadn}
{\cal V}_{ad}(\eta_n) = \{\,(u, q) \ : \  q\in F(\eta_n),\ u=u(q,\eta_n)\,\}.
\end{equation}
With these data we consider the following perturbation of  Problem ${\cal Q}$.

\medskip\noindent
{\bf Problem} ${\cal Q}_n$. {\it Given $\eta_n\in\Sigma$, find $(u^*_n, q^*_n)\in {\cal V}_{ad}(\eta_n)$ such that}
\begin{equation}\label{o1n}
{\cal L}(u^*_n,q^*_n)=\min_{(u,q)\in {\cal V}_{ad}(\eta_n)} {\cal L}(u,q).
\end{equation}


\medskip
We also consider the following auxiliary problem.

\medskip\noindent
{\bf Problem} ${\cal R}_n$. {\it Given $\eta_n\in\Sigma$, find $q^*_n\in F(\eta_n)$ such that}
\begin{equation}\label{o1pn}
J(q^*_n,\eta_n)=\min_{q\in F(\eta_n)} J(q,\eta_n).
\end{equation}

The proof of  Lemma \ref{lb} shows that Problem ${\cal R}_n$ has at least one solution $q^*_n$,  for each $n\in\mathbb{N}$.
Moreover, using  (\ref{Jmn}) and definition (\ref{vadn}) it is easy to see that
\begin{equation}\label{b2}
\left\{ \begin{array}{l}
(u^*_n,q^*_n) \ \mbox{is a solution of Problem}\ {\cal Q}_n\ \mbox{if and only if}  \\ [2mm]
q^*_n \ \mbox{is a solution of Problem}\ {\cal R}_n\ \mbox{and}\ u^*_n=u(q^*_n,\eta_n).
\end{array}\right.  
\end{equation}	
This implies that for each $n\in\mathbb{N}$, Problem ${\cal Q}_n$ has at least one solution $(u^*_n,q^*_n)$.

\medskip

In order to study the link between the solutions to Problems ${\cal Q}_n$ and Problem ${\cal Q}$ we consider the following assumptions.
\begin{eqnarray}
&&\label{co1}\quad \ \eta_n\rightharpoonup\eta\quad {\rm in}\quad \Theta.\\[3mm]
&&\label{co3}
\quad\ F(\eta_n)\stackrel{M}{\longrightarrow} F(\eta)\quad {\rm in}\quad Q.\\[3mm]
&&\label{co4}\left\{\begin{array}{l}
\mbox{For all sequences }\{u_k\}\subset X\mbox{ and
}\{q_k\}\subset
U \mbox{ such that}\\[2mm]
u_k\rightarrow u\mbox{\ \ in\ \ }X,\ \ q_{k}\rightharpoonup
q \mbox{\ \ in\ \ }
Q,\ \mbox{we have} \\[2mm]
\displaystyle\lim_{k\to
	\infty}\,{\cal L}(u_k,q_k)={\cal L}(u,q).
\end{array}\right.
\end{eqnarray}

Then, we have the following convergence result.

\begin{Theorem}\label{t5}  Assume  $(\ref{K})$--$(\ref{f})$, for any $p\in\Lambda$. Moreover, assume $(\ref{pi})$,  $(\ref{es})$, $(\ref{o2})$, $(\ref{o3})$  and, in addition,	 assume that either 
	 $(\ref{o4})$ or $(\ref{o5})$ holds. For each $n\in\mathbb{N}$,	
	 let $(u^*_n,q^*_n)$ be a solution of Problem ${\cal Q}_n$. Then, if $(\ref{co1})$--$(\ref{co4})$ hold,
	there exists a subsequence of the sequence $\{(u^*_n,q^*_n)\}$, again denoted by $\{(u^*_n,q^*_n)\}$, and an element $(u^*,q^*)\in X\times Q$, such that
	\begin{eqnarray}
	&&\label{c11}
	u^*_n\to u^*\quad \mbox{\rm in}\quad X,  
	\\[2mm]
	&&\label{c12}
	q^*_n\rightharpoonup q^*\quad \mbox{\rm in}\quad Q.  
	\end{eqnarray}
	Moreover, $(u^*,q^*)$  is a solution of Problem  ${\cal Q}$.
	
\end{Theorem}

The proof of Theorem \ref{t5} will be carried out in several steps that we present in what follows. To this end, everywhere below we assume that
the hypotheses of Theorem \ref{t5} are satisfied. The first step of the proof concerns the function (\ref{Jmn}) and is as follows.

\begin{Lemma}\label{lc} {\rm (i) } For all sequence 
		$\{\eta_k\}\subset
		\Sigma$ such that $\eta_k\rightharpoonup
		\eta \mbox{ in }
		\Theta$  and for all  $q\in U$, the equality below holds: 
	\begin{equation}\label{j3}\displaystyle\lim_{k\to
		\infty}\,J(q,\eta_k)=J(q,\eta).
	\end{equation}
	 {\rm (ii) } For all sequences $\{q_k\}\subset U$ and
	 	$\{\eta_k\}\subset
	 	\Sigma$ such that
	 	$q_k\to q\mbox{ in }Q$,\,\,$\eta_k\rightharpoonup
	 	\eta \mbox{ in }
	 	\Theta$ 
	 		the equality below holds: 
	 \begin{equation}\label{j4}\displaystyle\lim_{k\to
	 	\infty}\,[J(q_k,\eta_k)-J(q,\eta_k)]=0.
	 \end{equation}
\end{Lemma}

\medskip\noindent{\bf Proof.}  {\rm (i) } Let
$\{\eta_k\}\subset\Sigma$ be a sequence such that $\eta_k\rightharpoonup
\eta \mbox{ in } \Theta$ and let $q\in U$.  We use the definition (\ref{Jmn}) to see that
\begin{equation}\label{z1}
J(q,\eta_k)-J(q,\eta)={\cal L}(u(q,\eta_k),q)-{\cal L}(u(q,\eta),q).
\end{equation}
On the other hand, using the convergence $u(q,\eta_k)\to u(q,\eta)$ in $X$, guaranteed by assumption (\ref{es}), we deduce from  (\ref{co4}) that
\begin{equation}\label{z2}
\lim_{k\to
	\infty}[{\cal L}(u(q,\eta_k),q)-{\cal L}(u(q,\eta),q)]=0.
\end{equation}
We now combine relations (\ref{z1}) and (\ref{z2}) to see that the equality (\ref{j3}) holds.

\medskip
 {\rm (ii) } Assume now that $\{q_k\}\subset U$ and
$\{\eta_k\}\subset
\Sigma$ are two sequences  such that $q_k\to q\mbox{ in } Q$ and $\eta_k\rightharpoonup
\eta \mbox{ in } \Theta$. We write
\begin{eqnarray}
	&&\label{mx2}J(q_k,\eta_k)-J(q,\eta_k)={\cal L}(u(q_k,\eta_k),q_k)-{\cal L}(u(q,\eta_k),q)\\[2mm]
	&&\quad= {\cal L}(u(q_k,\eta_k),q_k)-{\cal L}(u(q,\eta),q_k)\nonumber\\[2mm]
	&&\qquad+{\cal L}(u(q,\eta),q_k)-{\cal L}(u(q,\eta),q)\nonumber\\[2mm]
&&\qquad\quad +{\cal L}(u(q,\eta),q)-{\cal L}(u(q,\eta_k),q)	\nonumber
\end{eqnarray}
and, since $u(q_k,\eta_k)\to u(q,\eta)$, $u(q,\eta_k)\to u(q,\eta)$, both in $X$, assumptions (\ref{o3})(b) and (\ref{co4}) imply that
\begin{eqnarray*}\label{mx3}
&&{\cal L}(u(q_k,\eta_k),q_k)-{\cal L}(u(q,\eta),q_k)\to 0,\\[2mm]
&&{\cal L}(u(q,\eta),q_k)-{\cal L}(u(q,\eta),q)\to 0,\\[2mm]
&&{\cal L}(u(q,\eta),q)-{\cal L}(u(q,\eta_k),q)	\to 0.
\end{eqnarray*}
Wo now use these convergences in equality
(\ref{mx2}) to  deduce that  (\ref{j4}) holds, which concludes the proof
\hfill$\Box$

\medskip
In the next step we prove the following convergence result concerning Problems ${\cal R}_n$ and  ${\cal R}$. 

\begin{Lemma}\label{ld} 
For each $n\in\mathbb{N}$,	
	let $q^*_n$ be a solution of Problem ${\cal R}_n$.
	Then, 
	there exists a subsequence of the sequence $\{q^*_n\}$, again denoted by $\{q^*_n\}$, and an element $q^*\in Q$, such that $(\ref{c12})$ holds.
	Moreover, $q^*$  is a solution of Problem  ${\cal R}$.
\end{Lemma}

\noindent{\bf Proof.}
We claim that the sequence $\{q_n^*\}$ is bounded in $Q$. This claim is obviously satisfied if we assume that
 $(\ref{o5})$ holds.  Assume in what follows that  $(\ref{o4})$ holds. If $\{q_n^*\}$ is not bounded in $Q$, then we can find a subsequence of the sequence $\{q_n^*\}$, again denoted by $\{q_n^*\}$,  such that $\|q_n^*\|_Q \to \infty$. Therefore,   using definition (\ref{Jmn}) and condition (\ref{o4}) we deduce that
 \[
 J(q_n^*,\eta_n)={\cal L}(u(q_n^*,\eta_n),q^*_n)\ge h(q_n^*)\to\infty,
 \]
which implies that
\begin{equation}\label{11m}
J(q^*_n,\eta_n)\to \infty.
\end{equation}

Let  $s$ be a given element in $F(\eta)$ and note that assumption (\ref{co3}) and condition $(M_1)$ in Definition \ref{M}  impliy that there exists a sequence
$\{s_n\}$ such that
$s_n\in F(\eta_n)$ for each $n\in \mathbb{N}$  and 
\begin{equation}\label{xm}
s_n\to s\quad \mbox{\rm in}\quad Q.  
\end{equation}
Moreover, since $q_n^*$ is a solution of Problem ${\cal R}_n$ we have 
$J(q_n^*,\eta_n)\le J(s_n,\eta_n)$
and, therefore,
\begin{equation}\label{13m}
J(q_n^*,\eta_n)\le [J(s_n,\eta_n)- J(s,\eta_n)]+ [J(s,\eta_n)- J(s,\eta)]+J(s,\eta)\quad\forall\, n\in\mathbb{N}.
\end{equation}

On the other hand, the convergences (\ref{xm}) and (\ref{co1}) allow us to use  equality
(\ref{j4}) in Lemma \ref{lc}(ii)  to find that
$J(s_n,\eta_n)-J(s,\eta_n)\to 0$
and, in addition,  equality (\ref{j3}) in Lemma \ref{lc}(i) shows that $J(s,\eta_n)- J(s,\eta)\to 0$.
Thus, (\ref{13m}) implies that the sequence $\{J(q_n^*,\eta_n)\}$ is bounded, which contradicts
(\ref{11m}).	We conclude from above that the sequence $\{q_n^*\}$ is bounded in $Q$ and, therefore, there exists a subsequence of the sequence $\{q_n^*\}$, again denoted by $\{q_n^*\}$, and an element $q^*\in Q$, such that (\ref{c12}) holds.

We now prove that $q^*$  is a solution of  Problem  ${\cal R}$. To this end
we recall that $q_n^*\in F(\eta_n)$, for all $n\in\mathbb{N}$. Therefore, using (\ref{c12}) and condition $(M_2)$ in Definition \ref{M}, guaranteed by assumption (\ref{co3}), we deduce that $q^*\in F(\eta)$.
Next, we consider an arbitrary element $s\in F(\eta)$ and, using
condition $(M_1)$, we know that there exists a   sequence
$\{s_n\}$ such that
$s_n\in F(\eta_n)$ for each $n\in \mathbb{N}$  and (\ref{xm}) holds. Since $q_n^*$ is the solution to Problem 
${\cal R}_n$ we have
$J(q_n^*,\eta_n)\le J(s_n,\eta_n)$
which implies that
\begin{equation}\label{14m}
0\le\,[J(s,\eta_n)-J(q_n^*,\eta_n)]+[J(s_n,\eta_n)-J(s,\eta_n)].
\end{equation} 
We now use the convergences (\ref{co1}), (\ref{c12}), (\ref{xm}), Lemma \ref{la} and Lemma \ref{lc}(ii)
to see that
\begin{equation}\label{15m}
\limsup_{n\to\infty}\,[J(s,\eta_n)-J(q^*_n,\eta_n)]\le J(s,\eta)-J(q^*,\eta),
\end{equation} 
\begin{equation}\label{16m}
\lim_{n\to\infty}\,[J(s_n,\eta_n)-J(s,\eta_n)]=0,
\end{equation} 	
respectively. We now pass to the upper limit (\ref{14m}) and use  (\ref{15m}), (\ref{16m}) to deduce that  
\[0\le J(s,\eta)-J(q^*,\eta)\]
and, since $q^*\in F(\eta)$, we deduce that
$q^*$  is a solution of  Problem  ${\cal R}$, which concludes the proof.
\hfill$\Box$

\medskip We now have all the ingredients to provide the proof of Theorem \ref{t5}.
	
\medskip\noindent{\bf Proof.} We first use (\ref{b2}) to see that, for each $n\in\mathbb{N}$, $q_n^*$ is a solution to Problem ${\cal R}_n$. Then, we use Lemma \ref{ld} to see that there exists a subsequence of the sequence $\{q^*_n\}$, again denoted by $\{q^*_n\}$, and an element $q^*\in Q$, such that $(\ref{c12})$ holds. Moreover, $q^*$  is a solution of Problem  ${\cal R}$. Let $u^*=u(q^*,\eta)$. Then, equivalence
(\ref{b1}) shows that $(u^*,q^*)$ is a solution of Problem  ${\cal Q}$. On the other hand, (\ref{b2})  implies that $u^*_n=u(q^*_n,\eta_n)$ and the convergences  (\ref{c12}), (\ref{co1}) show that $p^*_n=(q_n^*,\eta_n)\rightharpoonup p^*=(q^*,\eta)$ in $Z$.  It follows now from assumption (\ref{es}) that (\ref{c11}) holds, which concludes the proof.
\hfill$\Box$

\section{An elastic contact problem}\label{s5}
\setcounter{equation}0

The abstract results  in Sections \ref{s3}--\ref{s4} are useful
in the study of various mathematical models which describe the equilibrium of elastic bodies in frictional contact with a foundation. In this section
we provide an example of such model and, to this end, we need some  notations and preliminaries. For details on the material below we refer the reaer
to \cite{HS,SofMat}, for instance.

Let $d\in\{1,2,3\}$. We denote by $\mathbb{S}^d$ the space of second order symmetric tensors on $\mathbb{R}^d$ or, equivalently, the space of symmetric
matrices of order $d$.
We  recall that
inner product and norm on
$\mathbb{R}^d$ and $\mathbb{S}^d$ are defined by
\begin{eqnarray*}
	&&\bu\cdot \bv=u_i v_i\ ,\qquad
	\|\bv\|=(\bv\cdot\bv)^{\frac{1}{2}}\qquad
	\forall \,\bu, \bv\in \mathbb{R}^d,\\[0mm]
	&&\bsigma\cdot \btau=\sigma_{ij}\tau_{ij}\ ,\qquad
	\|\btau\|=(\btau\cdot\btau)^{\frac{1}{2}} \qquad \,\forall\,
	\bsigma,\btau\in\mathbb{S}^d,
\end{eqnarray*}
where the indices $i$, $j$ run between $1$ and $d$ and,
unless stated otherwise, the summation convention over repeated
indices is used. The zero element of the spaces $\mathbb{R}^d$ and $\mathbb{S}^d$ will be denoted by $\bzero$. 

Let $\Omega\subset\mathbb{R}^d$ be a domain with smooth boundary $\Gamma$. The boundary $\Gamma$ is divided into three
measurable disjoint parts $\Gamma_1$, $\Gamma_2$ and $\Gamma_3$ such that ${ meas}\,(\Gamma_1)>0$.  
A generic point in $\Omega\cup\Gamma$ will be denoted by $\bx=(x_i)$.
We use the
standard notation for Sobolev and Lebesgue spaces associated to
$\Omega$ and $\Gamma$. In particular, we use the spaces  $L^2(\Omega)^d$, $L^2(\Gamma_2)^d$,
$L^2(\Gamma_3)$  and $H^1(\Omega)^d$, endowed with their canonical inner products and associated norms.
Moreover, for an element $\bv\in H^1(\Omega)^d$ we still  write $\bv$ for the trace of
$\bv$ to $\Gamma$. In addition, we consider the 
space
\begin{eqnarray*}
	&&V=\{\,\bv\in H^1(\Omega)^d:\  \bv =\bzero\ \ {\rm on\ \ }\Gamma_1\,\},
\end{eqnarray*}
which is a real Hilbert space
endowed with the canonical inner product 
\begin{equation}
	(\bu,\bv)_V= \int_{\Omega}
	\bvarepsilon(\bu)\cdot\bvarepsilon(\bv)\,dx
\end{equation}
and the associated norm
$\|\cdot\|_V$. Here and below $\bvarepsilon$ 
represents the deformation operator, i.e.,
\[
\bvarepsilon(\bu)=(\varepsilon_{ij}(\bu)),\quad
\varepsilon_{ij}(\bu)=\frac{1}{2}\,(u_{i,j}+u_{j,i}),
\]  
where an index that follows a comma will represent the
partial derivative with respect to the corresponding component of $\bx$, e.g.,\ $u_{i,j}=\frac{\partial u_i}{\partial j}$.
The completeness of the space $V$ follows from the assumption
${ meas}\,(\Gamma_1)>0$ which allows the use of Korn's inequality.
We denote by $\bzero_V$ the zero element of $V$ and we recall that, for an element $\bv\in V$,  the  normal and tangential components on $\Gamma$
are given by
$v_\nu=\bv\cdot\bnu$ and $\bv_\tau=\bv-v_\nu\bnu$, respectively.  Finally, $V^*$ represents the dual of $V$ and
$\langle\cdot,\cdot\rangle$ denotes the duality pairing
 between $V^*$ and $V$.
We also denote by $\|\gamma\|$ the norm of the trace operator $\gamma:V\to L^2(\Gamma_3)^d$ and we recall the inequality
\begin{equation}\label{trace}
\|\bv\|_{L^2(\Gamma)^d}\leq \|\gamma\|\|\bv\|_{V}\qquad \forall\,
	\bv\in V.
	\end{equation}
In addition, we use the space $Y=L^2(\Omega)^d\times L^2(\Gamma_3)^d$ equiped with the canonical product topology and the operator $\pi:V\to Y$
defined by	
\begin{equation}
\label{pim}\pi \bv=(\bv,\gamma\bv)\qquad\forall\, \bv\in V.
\end{equation}

The contact model we consider in this section is constructed by using a function ${\cal F}$ and a set $B$ which satisfy
the following conditions.

\begin{equation}
\left\{\begin{array}{ll} {\cal F}\colon 
\mathbb{S}^d\to \mathbb{S}^d\ \mbox{is such that} \\ [1mm]
{\rm (a)\  there\ exists}\ L_{\cal F}>0\ {\rm such\ that}\\
{}\qquad \|{\cal F}\bvarepsilon_1-{\cal F}\bvarepsilon_2\|
\le L_{\cal F} \|\bvarepsilon_1-\bvarepsilon_2\|\\
{}\qquad \quad\mbox{for all} \ \ \bvarepsilon_1,\bvarepsilon_2
\in \mathbb{S}^d; 
\\ [1mm]
{\rm (b)\  there\ exists}\ m_{\cal F}>0\ {\rm such\ that}\\
{}\qquad ({\cal F}\bvarepsilon_1-{\cal F}\bvarepsilon_2)
\cdot(\bvarepsilon_1-\bvarepsilon_2)\ge m_{\cal F}\,
\|\bvarepsilon_1-\bvarepsilon_2\|^2\quad\\
{}\qquad\quad \mbox{for all} \ \ \bvarepsilon_1,
\bvarepsilon_2 \in \mathbb{S}^d.\\ [1mm]
\end{array}\right.
\label{F}
\end{equation}

 \begin{eqnarray}
 &&\label{smm}m_{\cal F}>\|\gamma\|^2.\\ [2mm]
 &&\label{B} B \ \ \mbox{is a closed  convex subset of} \ \mathbb{S}^d\ \mbox{such that}\ \bzero\in B.
 \end{eqnarray}
Note that assumption (\ref{smm}) allows us to find a constant $\widetilde{m}_0\in\R$ such that
\begin{equation}\label{m0}
m_{\cal F}-\|\gamma\|^2>\widetilde{m}_0>0.
\end{equation}
We denote by $P_B:\mathbb{S}^d\to B$ the projection operator on $B$ and, for given $\rho\ge 0$, we consider the functions
 $k_\rho:\mathbb{R}\to\mathbb{R}$,  $j_\rho:\mathbb{R}\to\mathbb{R}$ defined by
 \begin{equation}\label{kr}k_\rho(r)=\left\{\begin{array}{ll}0\hspace{17mm}{\rm if}\quad r<0,\\ [2mm]
 \,r\hspace{16mm}{\rm if}\quad 0\le r<\rho,\\ [2mm]
 2\rho-r\qquad{\rm if}\quad \rho\le r<2\rho,\\ [2mm]
 r-2\rho\qquad{\rm if}\quad r\ge 2\rho,\end{array}\right.
 \end{equation}
 \begin{equation}\label{jr}
 j_\rho(r)=\int_0^rk_\rho(s)\,ds\qquad\forall\,r\in\R.
 \end{equation}
It is easy to see that the function $k_\rho$ is Lipschitz continuous, yet  not monotone. As a result the function $j_\rho$ is not convex.

Let $Z=\mathbb{R}^4\times L^2(\Omega)^d\times L^2(\Gamma_3)^d$ be the product space endowed with the 
the norm
\[\|p\|_{Z}=|\omega|+|\mu|+|\rho|+|g|+\|\fb_0\|_{L^2(\Omega)^d}+\|\fb_2\|_{L^2(\Gamma_3)^d}\] 
for all $p=(\omega,\mu,\rho,g,\fb_0,\fb_2)\in Z$
and let 
$\Lambda$ the subset of $Z$ defined by
\begin{equation}
\label{L}\Lambda=\{\,p=(\omega,\mu,\rho,g,\fb_0,\fb_2)\in Z\ :\ \omega,\ \mu, \rho, \ g\ge 0,\ \mu\|\gamma\|^2\le \widetilde{m}_0\,\}.
\end{equation}
Note that $Z$ is a weakly closed subset of $Z$ and, moreover, inequality $\widetilde{m}_0>0$  in (\ref{m0}) guarantees that $Z$ is not empty.
Next, for any $p=(\omega,\mu,\rho,g,\fb_0,\fb_2)\in\Lambda$ we define the set $K_p$,
the operator $A_p$, the functions $\varphi_p$, $j_p$ and the element $\fb_p$, by equalities
\begin{eqnarray}
\label{8b0}&&K_p=\{\,\bv\in V\ |\ v_\nu \le g\ \  \hbox{a.e. on}\
\Gamma_3\,\},\\ [3mm]
&&A_p\colon V \to V^*,\quad
\label{8b1}\langle A_p\bu,\bv\rangle =\int_{\Omega}\cF\bvarepsilon(\bu)\cdot\bvarepsilon(\bv)\,dx\\ [2mm]
&&\qquad\qquad\qquad+\omega\int_{\Omega}\big(\bvarepsilon(\bu)-P_B\bvarepsilon(\bu)\big)\cdot\bvarepsilon(\bv)\,dx,\nonumber\\[2mm]
&&\varphi_p\colon V\times V \to \real, 
\label{8b3}\quad
\varphi_p(\bu,\bv)=\int_{\Gamma_3}\mu\, u_\nu^+\,\|\bv_\tau\|\,da,\\[2mm]
&&j_p\colon V\to\real, \quad
\label{8b5}j_p(\bv)=\int_{\Gamma_3}j_\rho(v_\nu)\, da,\\[2mm]
&&\fb_p\in Y,\quad\label{8ef}(\fb_p, \pi \bv)_Y
=\int_{\Omega}\fb_0\cdot\bv\,dx +
\int_{\Gamma_3}\fb_2\cdot\gamma \bv\,da,
\end{eqnarray}
for all $\bu,\bv\in V$. Here and below $r^+$ represents the positive part of $r$, i.e., $r^+=max\,\{r,0\}$. With these notation, for $p\in\Lambda$ given,  we consider the following problem.

\medskip\noindent
{\bf Problem} ${\cal S}$.
	{\it Find a displacement field $\bu=\bu(p)$ such that}
	\begin{eqnarray}
	&&\label{hv}\bu\in K_p,\qquad
	\langle A_p\bu,\bv - \bu \rangle + \varphi_p (\bu, \bv) - \varphi_p (\bu, \bu) + j^0_p(\bu;\bv -\bu)
	\\[2mm]
	&&\qquad\ge (\fb_p, \pi\bv - \pi\bu )_Y 
	\ \ \mbox{for all} \ \ \bv \in V.\nonumber
	\end{eqnarray}

Following the arguments in \cite{SMBOOK}, it can be shown that Problem ${\cal S}$ represents the variational formulation of a mathematical model which describes the equilibrium of an elastic body in frictional contact with a foundation, under
the action of external forces. It is assumed that the foundation is made of a rigid material covered by a layer of deformable material.
The data ${\cal F}$, $\omega$ and $B$ are related to the constitutive law,  while $\fb_0$ and $\fb_2$ denote the density of body forces and applied tractions which act on the body and the surface $\Gamma_2$, respectively. In addition, $\mu$ represents the coefficient of friction and $\rho$ is a given stiffness coefficient. Finally, $g$ represents the thickness of the deformable layer. Note that here we restrict ourselves to the homogenuous case, for the convenience of the reader. Nevertheless, we mention that the mathematical tools in Sections \ref{s3} and \ref{s4} allow us to extend our results below to the nonhomogeneous case in which ${\cal F}$, $\omega$, $\mu$ depend on the spatial variable.

Contact models which lead to inequality problems of the form (\ref{hv}) have been considered in the books \cite{MOSBOOK,SMBOOK}.
There, the reader can find the classical formulation of the models, including the mechanical assumptions which lead to their construction,  as well as their variational analysis, based on arguments similar to that we briefly resume below in this section.

\medskip
Our first result in the study of Problem ${\cal S}$ is the following.

\begin{Theorem}\label{t6}  
	Assume $(\ref{F})$--$(\ref{B})$.
	Then, for each $p=(\omega,\mu,\rho,g,\fb_0,\fb_2)\in\Lambda$ there exists a unique solution $\bu=\bu(p)$ to the Problem ${\cal S}$.
	Moreover, if the sequence  $\{p_n\}$ with $p_n=(\omega_n,\mu_n,\rho_n,g_n,\fb_{0n},\fb_{2n})\in\Lambda$ is such that $p_n\rightharpoonup p$ in $Z$, then $\bu(p_n) \to \bu(p)$ in $V$.
\end{Theorem}

\noindent
{\bf Proof.}  Let $p=(\omega,\mu,\rho,g,\fb_0,\fb_2)\in\Lambda$. For the existence and uniqueness part we apply Theorem~\ref{t0} on  the space $X=V$. To this end,	we remark that, obviously, condition (\ref{K}) is satisfied. Moreover, the operator $A_p$ defined by \eqref{8b1} satisfies condition (\ref{A}).
	Indeed, for $\bu$, $\bv$, $\bw\in V$, by assumption (\ref{F})(a)  and the properties of the projection operator $P_B$, we have
	\[
	\langle A_p\bu-A_p\bv,\bw\rangle_{V^*\times V} 
	\leq (L_{\cal F}+2\omega)\|\bu-\bv\|_V\|\bw\|_V.
	\]
	This proves that
	$$\|A_p\bu-A_p\bv\|_{V^*}\le( L_{\cal F}+2\omega)\|\bu-\bv\|_V,$$
	for all $\bu$, $\bv\in V$, which implies that $A_p$ is Lipschitz
	continuous. On the other hand, using  assumption (\ref{F})(b) and the nonexpansivity of the projector operator yields
	\begin{equation}
	\label{Astrmo}\langle A_p\bu-A_p\bv,\bu-\bv\rangle_{V^*\times V}
	\geq m_{\cal F}\|\bu-\bv\|_V^2, \nonumber
	\end{equation}
	for all $\bu$, $\bv\in V$. This shows that condition (\ref{A})(b) is satisfied  with
	$m_p=m_{\mathcal F}$.
	Since 
	$A_p$ is Lipschitz continuous and monotone,  
	it follows that $A_p$ is pseudomonotone and, therefore, (\ref{A})(a) holds.

	Next, for $\varphi_p$ defined by \eqref{8b3}, we use the trace inequality (\ref{trace})  to see that  (\ref{fi}) holds with $\alpha_p = \mu\|\gamma \|^2$.

	On the other hand, since $j_\rho'(r)=k_\rho(r)$ for all $r\in\R$, it follows that $j_\rho$ is a $C^1$ function and, therefore,  is a locally Lipschitz function. Moreover, it is easy to see that
	$|k_\rho|\le 2\rho+|r|$ for $r\in\R$ and the function $r\mapsto k_\rho(r)+r$ is nondecreasing. We use these properties and equality $j_\rho^0(r;s)=k_ \rho(r)s$, valid for all $r,\ s\in\R$, to see that the function $j_\rho$ satisfies condition (\ref{j}) on $X=\R$ with $c_{0\rho}=2\rho$,  $c_{1\rho}=1$ and  $\beta_{\rho}=1$.
	Therefore, following the arguments in \cite[p.219]{SMBOOK} we deduce that
	the function $j_p$ given by (\ref{8b5}) satisfies (\ref{j}) with 
	$c_{0p} = 2\rho\sqrt{2\, meas (\Gamma_3)}\, 
	\| \gamma \|$, 
	$c_{1p} = \sqrt{2} \,  \| \gamma \|^2$ 
	and
	$\beta_p = \| \gamma \|^2$.
	
	It follows from above that $\alpha_p+\beta_p=(\mu+1)\|\gamma \|^2$ and using definition (\ref{L}) and inequality (\ref{m0}) we deduce that $\alpha_p+\beta_p<m_{\cal F}=m_p$ which shows that the smallness condition (\ref{sm}) is satisfied.
    We also note that  conditions (\ref{f}) and (\ref{pi}) are obviously satisfied.

	Therefore, we are in position to use Theorem~\ref{t0}. In this way we deduce that there exists a unique element $\bu\in V$ such that
	(\ref{hv}) holds, which concludes the existence and uniqueness part of the theorem.
	
	For the continuous depencence part we use Theorem \ref{t1}. Assume that \[p_n=(\omega_n,\mu_n,\rho_n,g_n,\fb_{0n},\fb_{2n})\rightharpoonup p=(\omega,\mu,\rho,g,\fb_0,\fb_2)\quad{\rm in}\quad Z\]
	which implies that
	\begin{eqnarray}
	&&\label{w1}\omega_n\to \omega,\quad \mu_n\to \mu,\quad \rho_n\to \rho,\quad g_n\to g\quad{\rm in}\ \quad\R,\\ [2mm]
	&&\label{w2}\fb_{0n}\rightharpoonup \fb_0\quad{\rm in}\ \quad L^2(\Omega)^d,\quad \fb_{02}\rightharpoonup \fb_2\quad{\rm in}\ \quad L^2(\Omega)^d.
	\end{eqnarray}
	Using the convergences (\ref{w1}) it is easy to see that conditions (\ref{cvK})--(\ref{cvfi}) are satisfied with
	$F_n= 2\,|\omega_n-\omega|$, $\delta_n=0$, 
	$c_n(r)=\mu_n\|\gamma\|^2\,r$. Moreover, the compactness of the trace and the convergence $\mu_n\to\mu$ imply that (\ref{cvfin}) holds, too.
	
	On the other hand, since the function $j_\rho$ is regular, Lemma 8 in \cite{SMBOOK} guarantees that 
	\begin{equation}\label{q1}
	j_p^0(\bu;\bv)=\int_{\Gamma_3}j_\rho^0(u_\nu;v_\nu)\,da=\int_{\Gamma_3}k_\rho(u_\nu)v_\nu\,da\qquad\forall\,\bu,\,\bv\in V, \ p\in\Lambda.
	\end{equation}
Moreover, an elementary argument shows that
\begin{equation}\label{q2}
	|k_{\rho_n}(r)-k_\rho(r)|\le2\,|\rho_n-\rho|\qquad\forall\, r\in \mathbb{R}.
\end{equation}
 Therefore, if $\bu_n\rightharpoonup \bu$ and $\bv_n\to v$ in $V$, using (\ref{q1}) and (\ref{q2}) we deduce that
\begin{eqnarray*}
&&j_{p_n}^0(\bu_n;\bv_n-\bu_n)=\int_{\Gamma_3}k_{\rho_n}(u_{n\nu})(v_{n\nu}-u_{n\nu})\,da\\ [2mm]
&&\quad \le 2\,|\rho_n-\rho|\int_{\Gamma_3}|v_{n\nu}-u_{n\nu}|\,da+\int_{\Gamma_3}k_{\rho}(u_{n\nu})(v_{n\nu}-u_{n\nu})\,da.
\end{eqnarray*}
Then, the convergence $\rho_n\to \rho$, the compactness of the trace and the properties of the function $k_\rho$  imply that
\[ \limsup j^0_{p_n}(\bu_n; \bv_n -\bu_n) \le \int_{\Gamma_3}k_{\rho}(u_{\nu})(v_{\nu}-u_{\nu})\,da=j^0_p(\bu; \bv - \bu)\]
which shows that condition (\ref{cvj}) holds. 

Next, by standard compactness arguments it follows that the operator (\ref{pim}) satisfies condition (\ref{cvpi}). Moreover, the convergences (\ref{w2}) show that (\ref{cvf}) holds, too.
Finally, we recall that	 $\alpha_{p_n} = \mu_n\|\gamma \|^2$,
$c_{0p_n} = 2\rho_n\sqrt{2\, meas (\Gamma_3)}\, 
\| \gamma \|$, 
$c_{1p_n} = \sqrt{2} \,  \| \gamma \|^2$, 
$\beta_{p_n} = \| \gamma \|^2$ and $m_{p_n}=m_{\cal F}$.
Therefore, using (\ref{L}) and (\ref{m0}) we deduce that condition (\ref{cvs})	  is satisfied with $m_0=m_{p_n}-\widetilde{m}_0-\| \gamma \|^2$ and, obviously, (\ref{cvb}) holds.

It follows from above that we are in position to use Theorem \ref{t1} in order to deduce  that $\bu(p_n) \to \bu(p)$ in $V$, which concludes the proof.
\hfill$\Box$

\medskip
We now associate  to Problem ${\cal S}$ an optimal control problem. To this end we note that, with a permutation of the factors, the space $Z=\mathbb{R}^4\times L^2(\Omega)^d\times L^2(\Gamma_3)^d$ can be written in the form
$Z =Q\times \Theta$, 
where $Q=\mathbb{R}\times L^2(\Gamma_3)^d$  
and $\Theta=\mathbb{R}^3\times L^2(\Omega)^d$. The norms on these spaces will be denoted by $\|\cdot\|_{Q}$ and $\|\cdot\|_{\Theta}$, respectively. 
We endow $Z$ with the norm
\[\|p\|_{Z}=\|q\|_{Q}+\|\eta\|_{\Theta}\qquad\forall\, p=(q_1,\eta)\in Z.\] A generic element of $Z$ will be denoted by $p=(q,\eta)$ where $q=(g,\fb_2)\in Q$ and $\eta=(\omega,\mu,\rho,\fb_0)\in\Theta$. 

Let $U\subset Q$, $\Sigma\subset\Theta$, $\Lambda\subset Z$ be the sets given by
\begin{eqnarray*}
&& U=\{\,q=(g,\fb_2)\in Q\ :\ 0\le g \le g_0,\ \|\fb_2\|_{L^2(\Gamma_3)^d}\le h_0\, \}, \\[2mm]
&& \Sigma=\{\,\eta=(\omega,\mu,\rho,\fb_0)\in \Theta\ :\ \omega,\ \mu, \rho\ge 0,\ \rho\le\rho_0,\ \mu\|\gamma\|^2\le \widetilde{m}_0\,\}, \\[2mm]
&&\Lambda=U\times \Sigma,
\end{eqnarray*}
where $g_0>0$, $h_0>0$ and $\rho_0>0$ are given. Moreover, for each $\eta=(\omega,\mu,\rho,\fb_0)\in\Sigma$, let $F(\eta)\subset U$ be the subset given by
\begin{equation}
 F(\eta)=\{\,q=(g,\fb_2)\in U\ :\ \rho\le g\le\rho_0\}.
\end{equation}
Then, under the assumptions $(\ref{F})$--$(\ref{B})$, it follows from Theorem \ref{t6} that Problem  ${\cal S}$ has a uniques solution
$\bu=\bu(p)=\bu(q,\eta)$.
We define
the set of admissible pairs for  Problem ${\cal S}$
by equality
\begin{equation*}
{\cal V}_{ad}(\eta) = \{\,(\bu, q) \ : \  q\in F(\eta),\ \bu=\bu(q,\eta)\,\}.
\end{equation*}
Consider also
the cost functional ${\cal L}:V\times U\to\mathbb{R}$
given by

\begin{equation}
{\cal L}(\bu, q)=\int_{\Gamma_3}(u_\nu-\phi)^2\,da\qquad\forall\, \bu\in V,\quad q\in U,
\end{equation}
where $\phi\in L^2(\Gamma_3)$ is given.
Then, the optimal control problem we are interested in is the following.

\medskip\noindent
{\bf Problem} ${\cal T}$. {\it Given $\eta=(\omega,\mu,\rho,\fb_0)\in\Sigma$, find $(\bu^*, q^*)\in {\cal V}_{ad}(\eta)$ such that}
\begin{equation*}
{\cal L}(\bu^*,q^*)=\min_{(u,q)\in {\cal V}_{ad}(\eta)} {\cal L}(u,q).
\end{equation*}

With this choice, the mechanical interpretation of  Problem ${\cal T}$ is the following:   
given a contact process described by the variational-hemivariational inequality  (\ref{hv}) with the  data $q=(g,\fb_2)\in U$ and $\eta=(\omega,\mu,\rho,\fb_0)\in\Sigma$,
we are looking for a thickness $g\in[\rho,\rho_0] $ such that 
the normal component of the corresponding solution  is as close as possible, on  $\Gamma_3$, to the ``desired  normal displacement" $\phi$.

Note that, in this case assumptions (\ref{o2}), (\ref{o3})
and (\ref{o5}) are satisfied.  Therefore  Theorem \ref{t4} guarantees the existence of the solutions of the  optimal control  problem ${\cal T}$. Moreover, note that assumption (\ref{co1}) implies  the convergence (\ref{co3}) and, in addition,   (\ref{co4}) holds. Therefore, the convergence result stated in Theorem \ref{t5} can be applied in the study of Problem ${\cal T}$.

\section*{Acknowledgements}

\indent This research was supported by the National Natural Science Foundation of China (11771067) and the Applied Basic Project of Sichuan Province (2016JY0170).

\end{document}